\documentclass[11pt]{article}
\long\def\ig#1{\relax}
\ig{Thanks to Roberto Minio for this def'n.  Compare the def'n of
\comment in AMSTeX.}

\newcount \coefa
\newcount \coefb
\newcount \coefc
\newcount\tempcounta
\newcount\tempcountb
\newcount\tempcountc
\newcount\tempcountd
\newcount\xext
\newcount\yext
\newcount\xoff
\newcount\yoff
\newcount\gap%
\newcount\arrowtypea
\newcount\arrowtypeb
\newcount\arrowtypec
\newcount\arrowtyped
\newcount\arrowtypee
\newcount\height
\newcount\width
\newcount\xpos
\newcount\ypos
\newcount\run
\newcount\rise
\newcount\arrowlength
\newcount\halflength
\newcount\arrowtype
\newdimen\tempdimen
\newdimen\xlen
\newdimen\ylen
\newsavebox{\tempboxa}%
\newsavebox{\tempboxb}%
\newsavebox{\tempboxc}%

\makeatletter
\def\settypes(#1,#2,#3){\arrowtypea#1 \arrowtypeb#2 \arrowtypec#3}
\def\settoheight#1#2{\setbox\@tempboxa\hbox{#2}#1\ht\@tempboxa\relax}%
\def\settodepth#1#2{\setbox\@tempboxa\hbox{#2}#1\dp\@tempboxa\relax}%
\def\settokens[#1`#2`#3`#4]{%
     \def\tokena{#1}\def\tokenb{#2}\def\tokenc{#3}\def\tokend{#4}}
\def\setsqparms[#1`#2`#3`#4;#5`#6]{%
\arrowtypea #1
\arrowtypeb #2
\arrowtypec #3
\arrowtyped #4
\width #5
\height #6
}
\def\setpos(#1,#2){\xpos=#1 \ypos#2}

\def\bfig{\begin{picture}(\xext,\yext)(\xoff,\yoff)}
\def\efig{\end{picture}}

\def\putbox(#1,#2)#3{\put(#1,#2){\makebox(0,0){$#3$}}}

\def\settriparms[#1`#2`#3;#4]{\settripairparms[#1`#2`#3`1`1;#4]}%

\def\settripairparms[#1`#2`#3`#4`#5;#6]{%
\arrowtypea #1
\arrowtypeb #2
\arrowtypec #3
\arrowtyped #4
\arrowtypee #5
\width #6
\height #6
}

\def\resetparms{\settripairparms[1`1`1`1`1;500]\width 500}

\resetparms

\def\mvector(#1,#2)#3{
\put(0,0){\vector(#1,#2){#3}}%
\put(0,0){\vector(#1,#2){30}}%
}
\def\evector(#1,#2)#3{{
\arrowlength #3
\put(0,0){\vector(#1,#2){\arrowlength}}%
\advance \arrowlength by-30
\put(0,0){\vector(#1,#2){\arrowlength}}%
}}

\def\horsize#1#2{%
\settowidth{\tempdimen}{$#2$}%
#1=\tempdimen
\divide #1 by\unitlength
}

\def\vertsize#1#2{%
\settoheight{\tempdimen}{$#2$}%
#1=\tempdimen
\settodepth{\tempdimen}{$#2$}%
\advance #1 by\tempdimen
\divide #1 by\unitlength
}

\def\vertadjust[#1`#2`#3]{%
\vertsize{\tempcounta}{#1}%
\vertsize{\tempcountb}{#2}%
\ifnum \tempcounta<\tempcountb \tempcounta=\tempcountb \fi
\divide\tempcounta by2
\vertsize{\tempcountb}{#3}%
\ifnum \tempcountb>0 \advance \tempcountb by20 \fi
\ifnum \tempcounta<\tempcountb \tempcounta=\tempcountb \fi
}

\def\horadjust[#1`#2`#3]{%
\horsize{\tempcounta}{#1}%
\horsize{\tempcountb}{#2}%
\ifnum \tempcounta<\tempcountb \tempcounta=\tempcountb \fi
\divide\tempcounta by20
\horsize{\tempcountb}{#3}%
\ifnum \tempcountb>0 \advance \tempcountb by60 \fi
\ifnum \tempcounta<\tempcountb \tempcounta=\tempcountb \fi
}

\ig{ In this procedure, #1 is the paramater that sticks out all the way,
#2 sticks out the least and #3 is a label sticking out half way.  #4 is
the amount of the offset.}

\def\sladjust[#1`#2`#3]#4{%
\tempcountc=#4
\horsize{\tempcounta}{#1}%
\divide \tempcounta by2
\horsize{\tempcountb}{#2}%
\divide \tempcountb by2
\advance \tempcountb by-\tempcountc
\ifnum \tempcounta<\tempcountb \tempcounta=\tempcountb\fi
\divide \tempcountc by2
\horsize{\tempcountb}{#3}%
\advance \tempcountb by-\tempcountc
\ifnum \tempcountb>0 \advance \tempcountb by80\fi
\ifnum \tempcounta<\tempcountb \tempcounta=\tempcountb\fi
\advance\tempcounta by20
}

\def\putvector(#1,#2)(#3,#4)#5#6{{%
\xpos=#1
\ypos=#2
\run=#3
\rise=#4
\arrowlength=#5
\arrowtype=#6
\ifnum \arrowtype<0
    \ifnum \run=0
        \advance \ypos by-\arrowlength
    \else
        \tempcounta \arrowlength
        \multiply \tempcounta by\rise
        \divide \tempcounta by\run
        \ifnum\run>0
            \advance \xpos by\arrowlength
            \advance \ypos by\tempcounta
        \else
            \advance \xpos by-\arrowlength
            \advance \ypos by-\tempcounta
        \fi
    \fi
    \multiply \arrowtype by-1
    \multiply \rise by-1
    \multiply \run by-1
\fi
\ifnum \arrowtype=1
    \put(\xpos,\ypos){\vector(\run,\rise){\arrowlength}}%
\else\ifnum \arrowtype=2
    \put(\xpos,\ypos){\mvector(\run,\rise)\arrowlength}%
\else\ifnum\arrowtype=3
    \put(\xpos,\ypos){\evector(\run,\rise){\arrowlength}}%
\fi\fi\fi
}}

\def\putsplitvector(#1,#2)#3#4{
\xpos #1
\ypos #2
\arrowtype #4
\halflength #3
\arrowlength #3
\gap 140
\advance \halflength by-\gap
\divide \halflength by2
\ifnum \arrowtype=1
    \put(\xpos,\ypos){\line(0,-1){\halflength}}%
    \advance\ypos by-\halflength
    \advance\ypos by-\gap
    \put(\xpos,\ypos){\vector(0,-1){\halflength}}%
\else\ifnum \arrowtype=2
    \put(\xpos,\ypos){\line(0,-1)\halflength}%
    \put(\xpos,\ypos){\vector(0,-1)3}%
    \advance\ypos by-\halflength
    \advance\ypos by-\gap
    \put(\xpos,\ypos){\vector(0,-1){\halflength}}%
\else\ifnum\arrowtype=3
    \put(\xpos,\ypos){\line(0,-1)\halflength}%
    \advance\ypos by-\halflength
    \advance\ypos by-\gap
    \put(\xpos,\ypos){\evector(0,-1){\halflength}}%
\else\ifnum \arrowtype=-1
    \advance \ypos by-\arrowlength
    \put(\xpos,\ypos){\line(0,1){\halflength}}%
    \advance\ypos by\halflength
    \advance\ypos by\gap
    \put(\xpos,\ypos){\vector(0,1){\halflength}}%
\else\ifnum \arrowtype=-2
    \advance \ypos by-\arrowlength
    \put(\xpos,\ypos){\line(0,1)\halflength}%
    \put(\xpos,\ypos){\vector(0,1)3}%
    \advance\ypos by\halflength
    \advance\ypos by\gap
    \put(\xpos,\ypos){\vector(0,1){\halflength}}%
\else\ifnum\arrowtype=-3
    \advance \ypos by-\arrowlength
    \put(\xpos,\ypos){\line(0,1)\halflength}%
    \advance\ypos by\halflength
    \advance\ypos by\gap
    \put(\xpos,\ypos){\evector(0,1){\halflength}}%
\fi\fi\fi\fi\fi\fi
}

\def\putmorphism(#1)(#2,#3)[#4`#5`#6]#7#8#9{{%
\run #2
\rise #3
\ifnum\rise=0
  \puthmorphism(#1)[#4`#5`#6]{#7}{#8}{#9}%
\else\ifnum\run=0
  \putvmorphism(#1)[#4`#5`#6]{#7}{#8}{#9}%
\else
\setpos(#1)%
\arrowlength #7
\arrowtype #8
\ifnum\run=0
\else\ifnum\rise=0
\else
\ifnum\run>0
    \coefa=1
\else
   \coefa=-1
\fi
\ifnum\arrowtype>0
   \coefb=0
   \coefc=-1
\else
   \coefb=\coefa
   \coefc=1
   \arrowtype=-\arrowtype
\fi
\width=2
\multiply \width by\run
\divide \width by\rise
\ifnum \width<0  \width=-\width\fi
\advance\width by60
\if l#9 \width=-\width\fi
\putbox(\xpos,\ypos){#4}
{\multiply \coefa by\arrowlength
\advance\xpos by\coefa
\multiply \coefa by\rise
\divide \coefa by\run
\advance \ypos by\coefa
\putbox(\xpos,\ypos){#5} }%
{\multiply \coefa by\arrowlength
\divide \coefa by2
\advance \xpos by\coefa
\advance \xpos by\width
\multiply \coefa by\rise
\divide \coefa by\run
\advance \ypos by\coefa
\if l#9%
   \put(\xpos,\ypos){\makebox(0,0)[r]{$#6$}}%
\else\if r#9%
   \put(\xpos,\ypos){\makebox(0,0)[l]{$#6$}}%
\fi\fi }%
{\multiply \rise by-\coefc
\multiply \run by-\coefc
\multiply \coefb by\arrowlength
\advance \xpos by\coefb
\multiply \coefb by\rise
\divide \coefb by\run
\advance \ypos by\coefb
\multiply \coefc by70
\advance \ypos by\coefc
\multiply \coefc by\run
\divide \coefc by\rise
\advance \xpos by\coefc
\multiply \coefa by140
\multiply \coefa by\run
\divide \coefa by\rise
\advance \arrowlength by\coefa
\ifnum \arrowtype=1
   \put(\xpos,\ypos){\vector(\run,\rise){\arrowlength}}%
\else\ifnum\arrowtype=2
   \put(\xpos,\ypos){\mvector(\run,\rise){\arrowlength}}%
\else\ifnum\arrowtype=3
   \put(\xpos,\ypos){\evector(\run,\rise){\arrowlength}}%
\fi\fi\fi}\fi\fi\fi\fi}}

\def\puthmorphism(#1,#2)[#3`#4`#5]#6#7#8{{%
\xpos #1
\ypos #2
\width #6
\arrowlength #6
\putbox(\xpos,\ypos){#3\vphantom{#4}}%
{\advance \xpos by\arrowlength
\putbox(\xpos,\ypos){\vphantom{#3}#4}}%
\horsize{\tempcounta}{#3}%
\horsize{\tempcountb}{#4}%
\divide \tempcounta by2
\divide \tempcountb by2
\advance \tempcounta by30
\advance \tempcountb by30
\advance \xpos by\tempcounta
\advance \arrowlength by-\tempcounta
\advance \arrowlength by-\tempcountb
\putvector(\xpos,\ypos)(1,0){\arrowlength}{#7}%
\divide \arrowlength by2
\advance \xpos by\arrowlength
\vertsize{\tempcounta}{#5}%
\divide\tempcounta by2
\advance \tempcounta by20
\if a#8 %
   \advance \ypos by\tempcounta
   \putbox(\xpos,\ypos){#5}%
\else
   \advance \ypos by-\tempcounta
   \putbox(\xpos,\ypos){#5}%
\fi}}

\def\putvmorphism(#1,#2)[#3`#4`#5]#6#7#8{{%
\xpos #1
\ypos #2
\arrowlength #6
\arrowtype #7
\settowidth{\xlen}{$#5$}%
\putbox(\xpos,\ypos){#3}%
{\advance \ypos by-\arrowlength
\putbox(\xpos,\ypos){#4}}%
{\advance\arrowlength by-140
\advance \ypos by-70
\ifdim\xlen>0pt
   \if m#8%
      \putsplitvector(\xpos,\ypos){\arrowlength}{\arrowtype}%
   \else
      \putvector(\xpos,\ypos)(0,-1){\arrowlength}{\arrowtype}%
   \fi
\else
   \putvector(\xpos,\ypos)(0,-1){\arrowlength}{\arrowtype}%
\fi}%
\ifdim\xlen>0pt
   \divide \arrowlength by2
   \advance\ypos by-\arrowlength
   \if l#8%
      \advance \xpos by-40
      \put(\xpos,\ypos){\makebox(0,0)[r]{$#5$}}%
   \else\if r#8%
      \advance \xpos by40
      \put(\xpos,\ypos){\makebox(0,0)[l]{$#5$}}%
   \else
      \putbox(\xpos,\ypos){#5}%
   \fi\fi
\fi
}}

\def\topadjust[#1`#2`#3]{%
\yoff=10
\vertadjust[#1`#2`{#3}]%
\advance \yext by\tempcounta
\advance \yext by 10
}
\def\botadjust[#1`#2`#3]{%
\vertadjust[#1`#2`{#3}]%
\advance \yext by\tempcounta
\advance \yoff by-\tempcounta
}
\def\leftadjust[#1`#2`#3]{%
\xoff=0
\horadjust[#1`#2`{#3}]%
\advance \xext by\tempcounta
\advance \xoff by-\tempcounta
}
\def\rightadjust[#1`#2`#3]{%
\horadjust[#1`#2`{#3}]%
\advance \xext by\tempcounta
}
\def\rightsladjust[#1`#2`#3]{%
\sladjust[#1`#2`{#3}]{\width}%
\advance \xext by\tempcounta
}
\def\leftsladjust[#1`#2`#3]{%
\xoff=0
\sladjust[#1`#2`{#3}]{\width}%
\advance \xext by\tempcounta
\advance \xoff by-\tempcounta
}
\def\adjust[#1`#2;#3`#4;#5`#6;#7`#8]{%
\topadjust[#1``{#2}]
\leftadjust[#3``{#4}]
\rightadjust[#5``{#6}]
\botadjust[#7``{#8}]}

\def\putsquarep<#1>(#2)[#3;#4`#5`#6`#7]{{%
\setsqparms[#1]%
\setpos(#2)%
\settokens[#3]%
\puthmorphism(\xpos,\ypos)[\tokenc`\tokend`{#7}]{\width}{\arrowtyped}b%
\advance\ypos by \height
\puthmorphism(\xpos,\ypos)[\tokena`\tokenb`{#4}]{\width}{\arrowtypea}a%
\putvmorphism(\xpos,\ypos)[``{#5}]{\height}{\arrowtypeb}l%
\advance\xpos by \width
\putvmorphism(\xpos,\ypos)[``{#6}]{\height}{\arrowtypec}r%
}}

\def\putsquare{\@ifnextchar <{\putsquarep}{\putsquarep%
   <\arrowtypea`\arrowtypeb`\arrowtypec`\arrowtyped;\width`\height>}}
\def\square{\@ifnextchar< {\squarep}{\squarep
   <\arrowtypea`\arrowtypeb`\arrowtypec`\arrowtyped;\width`\height>}}
\def\squarep<#1>[#2`#3`#4`#5;#6`#7`#8`#9]{{
\setsqparms[#1]
\xext=\width                                          
\yext=\height                                         
\topadjust[#2`#3`{#6}]
\botadjust[#4`#5`{#9}]
\leftadjust[#2`#4`{#7}]
\rightadjust[#3`#5`{#8}]
\begin{picture}(\xext,\yext)(\xoff,\yoff)
\putsquarep<\arrowtypea`\arrowtypeb`\arrowtypec`\arrowtyped;\width`\height>%
(0,0)[#2`#3`#4`#5;#6`#7`#8`{#9}]%
\end{picture}%
}}

\def\putptrianglep<#1>(#2,#3)[#4`#5`#6;#7`#8`#9]{{%
\settriparms[#1]%
\xpos=#2 \ypos=#3
\advance\ypos by \height
\puthmorphism(\xpos,\ypos)[#4`#5`{#7}]{\height}{\arrowtypea}a%
\putvmorphism(\xpos,\ypos)[`#6`{#8}]{\height}{\arrowtypeb}l%
\advance\xpos by\height
\putmorphism(\xpos,\ypos)(-1,-1)[``{#9}]{\height}{\arrowtypec}r%
}}

\def\putptriangle{\@ifnextchar <{\putptrianglep}{\putptrianglep
   <\arrowtypea`\arrowtypeb`\arrowtypec;\height>}}
\def\ptriangle{\@ifnextchar <{\ptrianglep}{\ptrianglep
   <\arrowtypea`\arrowtypeb`\arrowtypec;\height>}}

\def\ptrianglep<#1>[#2`#3`#4;#5`#6`#7]{{
\settriparms[#1]%
\width=\height                         
\xext=\width                           
\yext=\width                           
\topadjust[#2`#3`{#5}]
\botadjust[#3``]
\leftadjust[#2`#4`{#6}]
\rightsladjust[#3`#4`{#7}]
\begin{picture}(\xext,\yext)(\xoff,\yoff)
\putptrianglep<\arrowtypea`\arrowtypeb`\arrowtypec;\height>%
(0,0)[#2`#3`#4;#5`#6`{#7}]%
\end{picture}%
}}

\def\putqtrianglep<#1>(#2,#3)[#4`#5`#6;#7`#8`#9]{{%
\settriparms[#1]%
\xpos=#2 \ypos=#3
\advance\ypos by\height
\puthmorphism(\xpos,\ypos)[#4`#5`{#7}]{\height}{\arrowtypea}a%
\putmorphism(\xpos,\ypos)(1,-1)[``{#8}]{\height}{\arrowtypeb}l%
\advance\xpos by\height
\putvmorphism(\xpos,\ypos)[`#6`{#9}]{\height}{\arrowtypec}r%
}}

\def\putqtriangle{\@ifnextchar <{\putqtrianglep}{\putqtrianglep
   <\arrowtypea`\arrowtypeb`\arrowtypec;\height>}}
\def\qtriangle{\@ifnextchar <{\qtrianglep}{\qtrianglep
   <\arrowtypea`\arrowtypeb`\arrowtypec;\height>}}

\def\qtrianglep<#1>[#2`#3`#4;#5`#6`#7]{{
\settriparms[#1]
\width=\height                         
\xext=\width                           
\yext=\height                          
\topadjust[#2`#3`{#5}]
\botadjust[#4``]
\leftsladjust[#2`#4`{#6}]
\rightadjust[#3`#4`{#7}]
\begin{picture}(\xext,\yext)(\xoff,\yoff)
\putqtrianglep<\arrowtypea`\arrowtypeb`\arrowtypec;\height>%
(0,0)[#2`#3`#4;#5`#6`{#7}]%
\end{picture}%
}}

\def\putdtrianglep<#1>(#2,#3)[#4`#5`#6;#7`#8`#9]{{%
\settriparms[#1]%
\xpos=#2 \ypos=#3
\puthmorphism(\xpos,\ypos)[#5`#6`{#9}]{\height}{\arrowtypec}b%
\advance\xpos by \height \advance\ypos by\height
\putmorphism(\xpos,\ypos)(-1,-1)[``{#7}]{\height}{\arrowtypea}l%
\putvmorphism(\xpos,\ypos)[#4``{#8}]{\height}{\arrowtypeb}r%
}}

\def\putdtriangle{\@ifnextchar <{\putdtrianglep}{\putdtrianglep
   <\arrowtypea`\arrowtypeb`\arrowtypec;\height>}}
\def\dtriangle{\@ifnextchar <{\dtrianglep}{\dtrianglep
   <\arrowtypea`\arrowtypeb`\arrowtypec;\height>}}

\def\dtrianglep<#1>[#2`#3`#4;#5`#6`#7]{{
\settriparms[#1]
\width=\height                         
\xext=\width                           
\yext=\height                          
\topadjust[#2``]
\botadjust[#3`#4`{#7}]
\leftsladjust[#3`#2`{#5}]
\rightadjust[#2`#4`{#6}]
\begin{picture}(\xext,\yext)(\xoff,\yoff)
\putdtrianglep<\arrowtypea`\arrowtypeb`\arrowtypec;\height>%
(0,0)[#2`#3`#4;#5`#6`{#7}]%
\end{picture}%
}}

\def\putbtrianglep<#1>(#2,#3)[#4`#5`#6;#7`#8`#9]{{%
\settriparms[#1]%
\xpos=#2 \ypos=#3
\puthmorphism(\xpos,\ypos)[#5`#6`{#9}]{\height}{\arrowtypec}b%
\advance\ypos by\height
\putmorphism(\xpos,\ypos)(1,-1)[``{#8}]{\height}{\arrowtypeb}r%
\putvmorphism(\xpos,\ypos)[#4``{#7}]{\height}{\arrowtypea}l%
}}

\def\putbtriangle{\@ifnextchar <{\putbtrianglep}{\putbtrianglep
   <\arrowtypea`\arrowtypeb`\arrowtypec;\height>}}
\def\btriangle{\@ifnextchar <{\btrianglep}{\btrianglep
   <\arrowtypea`\arrowtypeb`\arrowtypec;\height>}}

\def\btrianglep<#1>[#2`#3`#4;#5`#6`#7]{{
\settriparms[#1]
\width=\height                         
\xext=\width                           
\yext=\height                          
\topadjust[#2``]
\botadjust[#3`#4`{#7}]
\leftadjust[#2`#3`{#5}]
\rightsladjust[#4`#2`{#6}]
\begin{picture}(\xext,\yext)(\xoff,\yoff)
\putbtrianglep<\arrowtypea`\arrowtypeb`\arrowtypec;\height>%
(0,0)[#2`#3`#4;#5`#6`{#7}]%
\end{picture}%
}}

\def\putAtrianglep<#1>(#2,#3)[#4`#5`#6;#7`#8`#9]{{%
\settriparms[#1]%
\xpos=#2 \ypos=#3
{\multiply \height by2
\puthmorphism(\xpos,\ypos)[#5`#6`{#9}]{\height}{\arrowtypec}b}%
\advance\xpos by\height \advance\ypos by\height
\putmorphism(\xpos,\ypos)(-1,-1)[#4``{#7}]{\height}{\arrowtypea}l%
\putmorphism(\xpos,\ypos)(1,-1)[``{#8}]{\height}{\arrowtypeb}r%
}}

\def\putAtriangle{\@ifnextchar <{\putAtrianglep}{\putAtrianglep
   <\arrowtypea`\arrowtypeb`\arrowtypec;\height>}}
\def\Atriangle{\@ifnextchar <{\Atrianglep}{\Atrianglep
   <\arrowtypea`\arrowtypeb`\arrowtypec;\height>}}

\def\Atrianglep<#1>[#2`#3`#4;#5`#6`#7]{{
\settriparms[#1]
\width=\height                         
\xext=\width                           
\yext=\height                          
\topadjust[#2``]
\botadjust[#3`#4`{#7}]
\multiply \xext by2 
\leftsladjust[#3`#2`{#5}]
\rightsladjust[#4`#2`{#6}]
\begin{picture}(\xext,\yext)(\xoff,\yoff)%
\putAtrianglep<\arrowtypea`\arrowtypeb`\arrowtypec;\height>%
(0,0)[#2`#3`#4;#5`#6`{#7}]%
\end{picture}%
}}

\def\putAtrianglepairp<#1>(#2)[#3;#4`#5`#6`#7`#8]{{
\settripairparms[#1]%
\setpos(#2)%
\settokens[#3]%
\puthmorphism(\xpos,\ypos)[\tokenb`\tokenc`{#7}]{\height}{\arrowtyped}b%
\advance\xpos by\height
\advance\ypos by\height
\putmorphism(\xpos,\ypos)(-1,-1)[\tokena``{#4}]{\height}{\arrowtypea}l%
\putvmorphism(\xpos,\ypos)[``{#5}]{\height}{\arrowtypeb}m%
\putmorphism(\xpos,\ypos)(1,-1)[``{#6}]{\height}{\arrowtypec}r%
}}

\def\putAtrianglepair{\@ifnextchar <{\putAtrianglepairp}{\putAtrianglepairp%
   <\arrowtypea`\arrowtypeb`\arrowtypec`\arrowtyped`\arrowtypee;\height>}}
\def\Atrianglepair{\@ifnextchar <{\Atrianglepairp}{\Atrianglepairp%
   <\arrowtypea`\arrowtypeb`\arrowtypec`\arrowtyped`\arrowtypee;\height>}}

\def\Atrianglepairp<#1>[#2;#3`#4`#5`#6`#7]{{%
\settripairparms[#1]%
\settokens[#2]%
\width=\height
\xext=\width
\yext=\height
\topadjust[\tokena``]%
\vertadjust[\tokenb`\tokenc`{#6}]
\tempcountd=\tempcounta                       
\vertadjust[\tokenc`\tokend`{#7}]
\ifnum\tempcounta<\tempcountd                 
\tempcounta=\tempcountd\fi                    
\advance \yext by\tempcounta                  
\advance \yoff by-\tempcounta                 %
\multiply \xext by2 
\leftsladjust[\tokenb`\tokena`{#3}]
\rightsladjust[\tokend`\tokena`{#5}]%
\begin{picture}(\xext,\yext)(\xoff,\yoff)%
\putAtrianglepairp
<\arrowtypea`\arrowtypeb`\arrowtypec`\arrowtyped`\arrowtypee;\height>%
(0,0)[#2;#3`#4`#5`#6`{#7}]%
\end{picture}%
}}

\def\putVtrianglep<#1>(#2,#3)[#4`#5`#6;#7`#8`#9]{{%
\settriparms[#1]%
\xpos=#2 \ypos=#3
\advance\ypos by\height
{\multiply\height by2
\puthmorphism(\xpos,\ypos)[#4`#5`{#7}]{\height}{\arrowtypea}a}%
\putmorphism(\xpos,\ypos)(1,-1)[`#6`{#8}]{\height}{\arrowtypeb}l%
\advance\xpos by\height
\advance\xpos by\height
\putmorphism(\xpos,\ypos)(-1,-1)[``{#9}]{\height}{\arrowtypec}r%
}}

\def\putVtriangle{\@ifnextchar <{\putVtrianglep}{\putVtrianglep
   <\arrowtypea`\arrowtypeb`\arrowtypec;\height>}}
\def\Vtriangle{\@ifnextchar <{\Vtrianglep}{\Vtrianglep
   <\arrowtypea`\arrowtypeb`\arrowtypec;\height>}}

\def\Vtrianglep<#1>[#2`#3`#4;#5`#6`#7]{{
\settriparms[#1]
\width=\height                         
\xext=\width                           
\yext=\height                          
\topadjust[#2`#3`{#5}]
\botadjust[#4``]
\multiply \xext by2 
\leftsladjust[#2`#3`{#6}]
\rightsladjust[#3`#4`{#7}]
\begin{picture}(\xext,\yext)(\xoff,\yoff)%
\putVtrianglep<\arrowtypea`\arrowtypeb`\arrowtypec;\height>%
(0,0)[#2`#3`#4;#5`#6`{#7}]%
\end{picture}%
}}

\def\putVtrianglepairp<#1>(#2)[#3;#4`#5`#6`#7`#8]{{
\settripairparms[#1]%
\setpos(#2)%
\settokens[#3]%
\advance\ypos by\height
\putmorphism(\xpos,\ypos)(1,-1)[`\tokend`{#6}]{\height}{\arrowtypec}l%
\puthmorphism(\xpos,\ypos)[\tokena`\tokenb`{#4}]{\height}{\arrowtypea}a%
\advance\xpos by\height
\putvmorphism(\xpos,\ypos)[``{#7}]{\height}{\arrowtyped}m%
\advance\xpos by\height
\putmorphism(\xpos,\ypos)(-1,-1)[``{#8}]{\height}{\arrowtypee}r%
}}

\def\putVtrianglepair{\@ifnextchar <{\putVtrianglepairp}{\putVtrianglepairp%
    <\arrowtypea`\arrowtypeb`\arrowtypec`\arrowtyped`\arrowtypee;\height>}}
\def\Vtrianglepair{\@ifnextchar <{\Vtrianglepairp}{\Vtrianglepairp%
    <\arrowtypea`\arrowtypeb`\arrowtypec`\arrowtyped`\arrowtypee;\height>}}

\def\Vtrianglepairp<#1>[#2;#3`#4`#5`#6`#7]{{%
\settripairparms[#1]%
\settokens[#2]
\xext=\height                  
\width=\height                 
\yext=\height                  
\vertadjust[\tokena`\tokenb`{#4}]
\tempcountd=\tempcounta        
\vertadjust[\tokenb`\tokenc`{#5}]
\ifnum\tempcounta<\tempcountd%
\tempcounta=\tempcountd\fi
\advance \yext by\tempcounta
\botadjust[\tokend``]%
\multiply \xext by2
\leftsladjust[\tokena`\tokend`{#6}]%
\rightsladjust[\tokenc`\tokend`{#7}]%
\begin{picture}(\xext,\yext)(\xoff,\yoff)%
\putVtrianglepairp
<\arrowtypea`\arrowtypeb`\arrowtypec`\arrowtyped`\arrowtypee;\height>%
(0,0)[#2;#3`#4`#5`#6`{#7}]%
\end{picture}%
}}

\def\putCtrianglep<#1>(#2,#3)[#4`#5`#6;#7`#8`#9]{{%
\settriparms[#1]%
\xpos=#2 \ypos=#3
\advance\ypos by\height
\putmorphism(\xpos,\ypos)(1,-1)[``{#9}]{\height}{\arrowtypec}l%
\advance\xpos by\height
\advance\ypos by\height
\putmorphism(\xpos,\ypos)(-1,-1)[#4`#5`{#7}]{\height}{\arrowtypea}l%
{\multiply\height by 2
\putvmorphism(\xpos,\ypos)[`#6`{#8}]{\height}{\arrowtypeb}r}%
}}

\def\putCtriangle{\@ifnextchar <{\putCtrianglep}{\putCtrianglep
    <\arrowtypea`\arrowtypeb`\arrowtypec;\height>}}
\def\Ctriangle{\@ifnextchar <{\Ctrianglep}{\Ctrianglep
    <\arrowtypea`\arrowtypeb`\arrowtypec;\height>}}

\def\Ctrianglep<#1>[#2`#3`#4;#5`#6`#7]{{
\settriparms[#1]
\width=\height                          
\xext=\width                            
\yext=\height                           
\multiply \yext by2 
\topadjust[#2``]
\botadjust[#4``]
\sladjust[#3`#2`{#5}]{\width}
\tempcountd=\tempcounta                 
\sladjust[#3`#4`{#7}]{\width}
\ifnum \tempcounta<\tempcountd          
\tempcounta=\tempcountd\fi              
\advance \xext by\tempcounta            
\advance \xoff by-\tempcounta           %
\rightadjust[#2`#4`{#6}]
\begin{picture}(\xext,\yext)(\xoff,\yoff)%
\putCtrianglep<\arrowtypea`\arrowtypeb`\arrowtypec;\height>%
(0,0)[#2`#3`#4;#5`#6`{#7}]%
\end{picture}%
}}

\def\putDtrianglep<#1>(#2,#3)[#4`#5`#6;#7`#8`#9]{{%
\settriparms[#1]%
\xpos=#2 \ypos=#3
\advance\xpos by\height \advance\ypos by\height
\putmorphism(\xpos,\ypos)(-1,-1)[``{#9}]{\height}{\arrowtypec}r%
\advance\xpos by-\height \advance\ypos by\height
\putmorphism(\xpos,\ypos)(1,-1)[`#5`{#8}]{\height}{\arrowtypeb}r%
{\multiply\height by 2
\putvmorphism(\xpos,\ypos)[#4`#6`{#7}]{\height}{\arrowtypea}l}%
}}

\def\putDtriangle{\@ifnextchar <{\putDtrianglep}{\putDtrianglep
    <\arrowtypea`\arrowtypeb`\arrowtypec;\height>}}
\def\Dtriangle{\@ifnextchar <{\Dtrianglep}{\Dtrianglep
   <\arrowtypea`\arrowtypeb`\arrowtypec;\height>}}

\def\Dtrianglep<#1>[#2`#3`#4;#5`#6`#7]{{
\settriparms[#1]
\width=\height                         
\xext=\height                          
\yext=\height                          
\multiply \yext by2 
\topadjust[#2``]
\botadjust[#4``]
\leftadjust[#2`#4`{#5}]
\sladjust[#3`#2`{#5}]{\height}
\tempcountd=\tempcountd                
\sladjust[#3`#4`{#7}]{\height}
\ifnum \tempcounta<\tempcountd         
\tempcounta=\tempcountd\fi             
\advance \xext by\tempcounta           %
\begin{picture}(\xext,\yext)(\xoff,\yoff)
\putDtrianglep<\arrowtypea`\arrowtypeb`\arrowtypec;\height>%
(0,0)[#2`#3`#4;#5`#6`{#7}]%
\end{picture}%
}}

\def\setrecparms[#1`#2]{\width=#1 \height=#2}%
\def\recursep<#1`#2>[#3;#4`#5`#6`#7`#8]{{%
\width=#1 \height=#2
\settokens[#3]
\settowidth{\tempdimen}{$\tokena$}
\ifdim\tempdimen=0pt
  \savebox{\tempboxa}{\hbox{$\tokenb$}}%
  \savebox{\tempboxb}{\hbox{$\tokend$}}%
  \savebox{\tempboxc}{\hbox{$#6$}}%
\else
  \savebox{\tempboxa}{\hbox{$\hbox{$\tokena$}\times\hbox{$\tokenb$}$}}%
  \savebox{\tempboxb}{\hbox{$\hbox{$\tokena$}\times\hbox{$\tokend$}$}}%
  \savebox{\tempboxc}{\hbox{$\hbox{$\tokena$}\times\hbox{$#6$}$}}%
\fi
\ypos=\height
\divide\ypos by 2
\xpos=\ypos
\advance\xpos by \width
\xext=\xpos \yext=\height
\topadjust[#3`\usebox{\tempboxa}`{#4}]%
\botadjust[#5`\usebox{\tempboxb}`{#8}]%
\sladjust[\tokenc`\tokenb`{#5}]{\ypos}%
\tempcountd=\tempcounta
\sladjust[\tokenc`\tokend`{#5}]{\ypos}%
\ifnum \tempcounta<\tempcountd
\tempcounta=\tempcountd\fi
\advance \xext by\tempcounta
\advance \xoff by-\tempcounta
\rightadjust[\usebox{\tempboxa}`\usebox{\tempboxb}`\usebox{\tempboxc}]%
\bfig
\putCtrianglep<-1`1`1;\ypos>(0,0)[`\tokenc`;#5`#6`{#7}]%
\puthmorphism(\ypos,0)[\tokend`\usebox{\tempboxb}`{#8}]{\width}{-1}b%
\puthmorphism(\ypos,\height)[\tokenb`\usebox{\tempboxa}`{#4}]{\width}{-1}a%
\advance\ypos by \width
\putvmorphism(\ypos,\height)[``\usebox{\tempboxc}]{\height}1r%
\efig
}}

\def\recurse{\@ifnextchar <{\recursep}{\recursep<\width`\height>}}

\def\puttwohmorphisms(#1,#2)[#3`#4;#5`#6]#7#8#9{{%
\puthmorphism(#1,#2)[#3`#4`]{#7}0a
\ypos=#2
\advance\ypos by 20
\puthmorphism(#1,\ypos)[\phantom{#3}`\phantom{#4}`#5]{#7}{#8}a
\advance\ypos by -40
\puthmorphism(#1,\ypos)[\phantom{#3}`\phantom{#4}`#6]{#7}{#9}b
}}

\def\puttwovmorphisms(#1,#2)[#3`#4;#5`#6]#7#8#9{{%
\putvmorphism(#1,#2)[#3`#4`]{#7}0a
\xpos=#1
\advance\xpos by -20
\putvmorphism(\xpos,#2)[\phantom{#3}`\phantom{#4}`#5]{#7}{#8}l
\advance\xpos by 40
\putvmorphism(\xpos,#2)[\phantom{#3}`\phantom{#4}`#6]{#7}{#9}r
}}

\def\puthcoequalizer(#1)[#2`#3`#4;#5`#6`#7]#8#9{{%
\setpos(#1)%
\puttwohmorphisms(\xpos,\ypos)[#2`#3;#5`#6]{#8}11%
\advance\xpos by #8
\puthmorphism(\xpos,\ypos)[\phantom{#3}`#4`#7]{#8}1{#9}
}}

\def\putvcoequalizer(#1)[#2`#3`#4;#5`#6`#7]#8#9{{%
\setpos(#1)%
\puttwovmorphisms(\xpos,\ypos)[#2`#3;#5`#6]{#8}11%
\advance\ypos by -#8
\putvmorphism(\xpos,\ypos)[\phantom{#3}`#4`#7]{#8}1{#9}
}}

\def\putthreehmorphisms(#1)[#2`#3;#4`#5`#6]#7(#8)#9{{%
\setpos(#1) \settypes(#8)
\if a#9 %
     \vertsize{\tempcounta}{#5}%
     \vertsize{\tempcountb}{#6}%
     \ifnum \tempcounta<\tempcountb \tempcounta=\tempcountb \fi
\else
     \vertsize{\tempcounta}{#4}%
     \vertsize{\tempcountb}{#5}%
     \ifnum \tempcounta<\tempcountb \tempcounta=\tempcountb \fi
\fi
\advance \tempcounta by 60
\puthmorphism(\xpos,\ypos)[#2`#3`#5]{#7}{\arrowtypeb}{#9}
\advance\ypos by \tempcounta
\puthmorphism(\xpos,\ypos)[\phantom{#2}`\phantom{#3}`#4]{#7}{\arrowtypea}{#9}
\advance\ypos by -\tempcounta \advance\ypos by -\tempcounta
\puthmorphism(\xpos,\ypos)[\phantom{#2}`\phantom{#3}`#6]{#7}{\arrowtypec}{#9}
}}

\def\putarc(#1,#2)[#3`#4`#5]#6#7#8{{%
\xpos #1
\ypos #2
\width #6
\arrowlength #6
\putbox(\xpos,\ypos){#3\vphantom{#4}}%
{\advance \xpos by\arrowlength
\putbox(\xpos,\ypos){\vphantom{#3}#4}}%
\horsize{\tempcounta}{#3}%
\horsize{\tempcountb}{#4}%
\divide \tempcounta by2
\divide \tempcountb by2
\advance \tempcounta by30
\advance \tempcountb by30
\advance \xpos by\tempcounta
\advance \arrowlength by-\tempcounta
\advance \arrowlength by-\tempcountb
\halflength=\arrowlength \divide\halflength by 2
\divide\arrowlength by 5
\put(\xpos,\ypos){\bezier{\arrowlength}(0,0)(50,50)(\halflength,50)}
\ifnum #7=-1 \put(\xpos,\ypos){\vector(-3,-2)0} \fi
\advance\xpos by \halflength
\put(\xpos,\ypos){\xpos=\halflength \advance\xpos by -50
   \bezier{\arrowlength}(0,50)(\xpos,50)(\halflength,0)}
\ifnum #7=1 {\advance \xpos by
   \halflength \put(\xpos,\ypos){\vector(3,-2)0}} \fi
\advance\ypos by 50
\vertsize{\tempcounta}{#5}%
\divide\tempcounta by2
\advance \tempcounta by20
\if a#8 %
   \advance \ypos by\tempcounta
   \putbox(\xpos,\ypos){#5}%
\else
   \advance \ypos by-\tempcounta
   \putbox(\xpos,\ypos){#5}%
\fi
}}

\makeatother

\usepackage{makeidx}  

\makeatother

\usepackage{amsmath}
\usepackage{qtree}
\usepackage{dsfont}
\usepackage{stmaryrd}
\usepackage{enumerate}
\usepackage{hyperref}

\usepackage{lscape}

\def\stl{{{\bf st}^l}}
\def\str{{{\bf st}^r}}

\def\pu{{{\bf pile'up}}}
\def\pul{{{\bf pile'up}^l}}
\def\pur{{{\bf pile'up}^r}}

\def\cps{{{\bf CPS}}}
\def\cpsl{{{\bf CPS}^l}}
\def\cpsr{{{\bf CPS}^r}}

\def\eps{{{\bf eps}}}
\def\epsl{{{\bf eps}^l}}
\def\epsr{{{\bf eps}^r}}

\def\mos{{{\bf mos}}}
\def\mosl{{{\bf mos}^l}}
\def\mosr{{{\bf mos}^r}}

\def\Id{{{\bf Id}}}
\def\shl{{{\bf shuf}^l}}
\def\shr{{{\bf shuf}^r}}
\def\sh{{{\bf shuf}}}

\def\lk{\langle}
\def\rk{\rangle}

\newcommand{\lra}{\longrightarrow}

\newcommand{\ra}{\rightarrow}
\newcommand{\la}{\leftarrow}

\newcommand{\Ra}{\Rightarrow}

\def\la{{\lambda}}

\def\bt{{\bf t}}

\def\cC{{\cal C}}

\def\cF{{\cal F}}

\def\cP{{\cal P}}

\def\vect[#1]{\vec{#1}}

\begin{document}

\title{Continuation semantics for multi-quantifier sentences: operation-based approaches}
\author{Justyna Grudzi\'{n}ska and Marek Zawadowski}

\maketitle

\begin{abstract}
\noindent Classical scope-assignment strategies for multi-quantifier sentences involve quantifier phrase (QP)-movement (e.g., \cite{may}, \cite{may85}). More recent continuation-based approaches provide a compelling alternative, for they interpret QP's in situ - without resorting to Logical Forms or any structures beyond the overt syntax. The continuation-based strategies can be divided into two groups: those that locate the source of scope-ambiguity in the rules of semantic composition (e.g., \cite{barker}) and those that attribute it to the lexical entries for the quantifier words (e.g., \cite{kiselyov:shan}, \cite{barker:shan}). In this paper, we focus on the former operation-based approaches and the nature of the semantic operations involved. More specifically, we discuss three such possible operation-based strategies for multi-quantifier sentences, together with their relative merits and costs.


\end{abstract}

\section{Introduction}

Multi-quantifier sentences have been known to be ambiguous with different readings corresponding to how various quantifier phrases (QPs) are semantically related in the sentence. For example,
\begin{enumerate}[(1)]
  \item Some teacher gave every student most books
\end{enumerate}
admits of six different readings, and in general a simple sentence with $n$ QPs will be (at least) $n!$ ways ambiguous (we only consider readings where QPs are linearly ordered - what we will call asymmetric readings).

Sentence (1) can be represented as a Syntactic (Surface Structure) Tree and the challenge is to obtain corresponding Semantic (Computation) Trees that compute the truth-value of the sentence in each of its readings
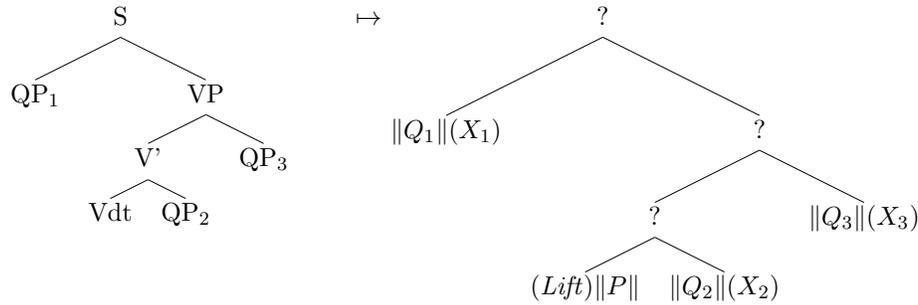
\begin{flushleft}
\scalebox{0.90}{
\Tree [.S [.QP$_1$ ][.VP [.V' [.Vdt ][.QP$_2$ ]][.QP$_3$ ]]]
\hskip 1 cm $\mapsto$
\Tree [.? [.$\|Q_1\|(X_1)$ ][.? [.? [.(\textit{Lift})$\|P\|$ ][.$\|Q_2\|(X_2)$ ]][.$\|Q_3\|(X_3)$ ]]]
\medskip
}
\end{flushleft}

\noindent We think here of Computation Trees by analogy with mathematical expressions, e.g.
\[ ((2-7)-8)+((12+5):7) \]
can be represented as
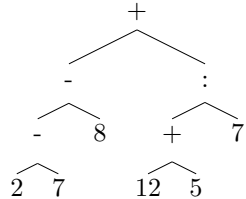
\scalebox{0.90}{
\par
\hskip 0,9 cm \Tree [.+ [.- [.- 2 7 ] 8 ][.: [.+ 12 5 ] 7 ]]
\medskip
}

\noindent i.e. a labeled binary tree where the leaves of this tree are labeled with values and the internal nodes are labeled with operations that will be applied in the computation to the values obtained from the computations of the left and right subtrees.

The process transforming Surface Structure Trees into the Semantic (Computation) Trees has been variously implemented in linguistics. Below we provide a list of desiderata relevant for the variety of the continuation-based semantics to be discussed in this paper.\\
\noindent (1) The semantics should be empirically adequate, i.e., it should allow us to calculate the truth-value of a given sentence in each of its readings.\\
\noindent (2) The semantics should be in situ, i.e., Computation Trees should have the same shape as Surface Structure Trees (with leaves labeled with interpretations of lexical items and inner nodes labeled with semantic operations).\\
\noindent (3) The semantic operations used should be kept as simple as possible.\\
\noindent (4) The semantic operations used should be kept as uniform as possible.\\
\noindent (5) The interpretation process should operate `on the fly', i.e., the particular reading(s) of a sentence should be determined as late as possible.

Classical movement analyses (involving quantifier phrase (QP)-movement, e.g., \cite{may}, \cite{may85}) generally meet desiderata (1), (3) and (4). More recent continuation-based approaches provide a compelling alternative, for they interpret QP's in situ - without resorting to Logical Forms or any structures beyond the overt syntax. They achieve this at a certain price though (including loss in simplicity/uniformity of the semantic operations used) - as summarized by the table below

\[
\begin{array}{|c|c|c|c|c|c|} \hline

  \textit{} & \textit{}  & \textit{} & \textit{} & &\\
  \textit{\textbf{Strategy}} & \textit{\textbf{A}} & \textit{\textbf{B}} & \textit{\textbf{C}}  & \textit{\textbf{D}} & \textit{\textbf{E}}\\
  \textit{} & \textit{Movement}  &\textit{Polyadic}  & \textit{Cont.-} & \textit{Minimal.} & \textit{Maximal.} \\
  \textit{} & \textit{Strategy} & \textit{Strategy}  & \textit{Based} & \textit{Augment.} & \textit{Uniformiz.} \\
   \textit{} & \textit{(May 77)}  & \textit{(May 85)}& \textit{Approach} & \textit{\textbf{Strat C}} & \textit{\textbf{Strat D}} \\
    \textit{} &   &   &\textit{(Barker 02)} & \textit{} & \textit{} \\ \hline
  \textit{} & \textit{} & & \textit{} & \textit{} & \\
  \textit{empirically} & \textit{YES}  & \textit{YES}  &\textit{NO} & \textit{YES} & \textit{YES}\\
  \textit{adequate} & \textit{}  & \textit{}  &\textit{} & \textit{} & \\
    \textit{} & \textit{} & \textit{} &\textit{} & \textit{ } & \\ \hline
   \textit{} & \textit{} & \textit{} & \textit{} & \textit{} & \\
  \textit{in situ} & \textit{NO} & \textit{NO} &\textit{YES}  & \textit{YES} & \textit{YES}\\
  \textit{} & \textit{} & \textit{} &\textit{} & \textit{ } & \\ \hline
   \textit{} & \textit{} & \textit{} &\textit{} & \textit{} & \textit{}\\
 \textit{semantic} & \textit{ $\mos$} & $\mos$ & \textit{$\cps$}  & \textit{$\cps$} & \textit{$\cps$} \\
  \textit{operations} & & $\pu$ & \textit{} & \textit{$\eps$} & $\eps$ \\
   \textit{} &  &  &\textit{} & $\sh$ & $\sh$\\
   \textit{} & \textit{} & \textit{} & \textit{} & \textit{} & \textit{(across}\\
    \textit{} & \textit{} & \textit{} &\textit{} & \textit{} &\textit{the board)}  \\
     \textit{their} & \textit{} & \textit{} &\textit{} & \textit{} &  \\
   \textit{simplicity} & \textit{YES} &\textit{?} & \textit{??}  & \textit{???} & \textit{???}\\
  \textit{\&} & \textit{} & \textit{} &\textit{} & \textit{ } & \\
    \textit{uniformity} & \textit{YES} &\textit{YES} & \textit{YES} & \textit{NO} & \textit{YES}\\
    \textit{} & \textit{} & \textit{} & \textit{} &\textit{ } & \\ \hline
    \textit{} & \textit{}  & \textit{} & \textit{} &\textit{} & \\
     \textit{`on the fly'} & \textit{NO} & \textit{NO} &\textit{YES} & \textit{?} & \textit{??}\\
     \textit{process} & \textit{} & \textit{} & \textit{} &\textit{ } & \\
    &  &  &  &  &\\ \hline

 \end{array}
\]
\vskip 4mm
\noindent In this paper, we first briefly recall the definitions of the semantic operations used in strategies \textbf{A}, \textbf{B} and \textbf{C}, i.e., $\mos$'es, $\pu$'es, $\cps$'es (for the details, see \cite{GZ}). We then discuss the three in situ strategies \textbf{C}, \textbf{D} and \textbf{E}, together with their relative merits and costs with respect to the desiderata introduced.

\section{Semantic Operations}

As noticed in \cite{barker} and \cite{groote}, a generalized quantifier on a set $X$ is an element of $\cC(X)$, the value of the continuation monad $\cC$ on $X$. Continuation-based semantics make heavy use of the computational machinery connected to the monad $\cC$ (its strength and derived operations). Here we only recall the definitions of the continuation monad, strengths and derived operations. For more details, we refer the reader to \cite{GZ}.

\subsection{Continuation Monad and Strengths}

We shall be working in the (cartesian closed) category of sets $Set$. The category $Set$ has sets as objects. A morphism in $Set$ from an object (set) X to an object (set) Y is a function $f: X\ra Y$ from $X$ to $Y$. For the general notion of a strong monad, see the appendix.
\vskip 4mm
\noindent {\bf Continuation monad $\cC$ - endofunctor}
 \begin{itemize}
   \item At the level of objects, it is just twice iterated power-set construction, i.e. for set $X$, $\cC(X)=\cP^2(X)$.\\
   ($\bt=\{ 0,1\}$, $\cP(X)=X\Ra \bt$ - powerset of $X$)
   \item At the level of morphisms, it is an inverse image of an inverse image, i.e., function $f:X\ra Y$ induces an inverse image function between powersets \[ \cP(f)=f^{-1}: \cP(Y)\ra \cP(X) \]
       \[ h\mapsto h\circ f, \hskip 2mm \cP(f)= \la h_{:\cP(Y)}.\la x_{:X}. h(f\, x) \]
  Taking again an inverse image function, we have
      \[ \cC(f)=\cP(f^{-1}) : \cC(X)=\cP^2(X)\ra \cP^2(Y)=\cC(Y) \]
  \[ Q \mapsto Q\circ f^{-1} , \hskip 2mm \cC(f)(Q)= \lambda h_{:\cP(Y)}.  Q(\lambda x_{:X}. h( f\, x))   \]
  for $Q\in \cC(X)$.
 \end{itemize}
 {\bf Continuation monad $\cC$ - natural transformations}
 \begin{itemize}
   \item  The {\bf unit}
 \[ \eta_X : X\ra \cC(X)\]
 is given by
  \[ \eta_X(x)= \lambda h_{:\cP(X)}. h(x)\]
  for $x\in X$.\\
  (lifts elements of $X$ as $\cC$-computations.)
   \item The {\bf multiplication} \[ \mu_X: \cC^{2}(X)\lra \cC(X)\] is given by
  \[\mu_X(\cF)(h)=\cF(\la D_{: \cC(X)}.D(h))\]
  for $\cF\in\cC^{2}(X)$ and  $h\in \cP(X)$.\\
  (flattens $\cC$-computations on $\cC$-computations to $\cC$-computations.)
 \end{itemize}
{\bf Continuation monad $\cC$ - strengths}\\
For the continuation monad, the {\bf left strength}  is
\[\stl: \cC(X)\times Y \lra \cC(X\times Y)\]
\[\stl(N,y)= \la c_{:\cP(X\times Y)}. N(\la x_{:X}. c(x,y))\]
for $N\in\cC(X)$ and $y\in Y$,\\
and the \textbf{right strength} is
\[\str: X\times \cC(Y) \lra \cC(X\times Y)\]
\[\str(x,M)=\la c_{: \cP(X\times Y)}. M(\lambda y_{:Y}.c(x,y))\]
for $x\in X$ and $M\in\cC(Y)$ .\\
Strengths allow  to lift pairs of $\cC$-computations to $\cC$-computations on products.

\subsection{Derived Operations}

{\bf $\pu$ operations}\\
Using both strengths, we can define $\pu$-operations.\\
For $M\in \cC(X)$ and $N\in \cC(Y)$, we have
\[\pul:\cC(X)\times \cC(Y)\ra \cC(X\times Y) \]
\[ \pul(M,N)=\la c_{:\cP(X\times Y)}. M(\la x_{:X}. N(\la y_{:Y} c(x,y))\]
and
\[\pur:\cC(X)\times \cC(Y)\ra \cC(X\times Y). \]
\[ \pur(M,N)= \la c_{:\cP(X\times Y)}. N(\la y_{:Y}. M(\la x_{:X} c(x,y)).\]
Thus $\pu$-operations put (interpretations of) quantifiers in order, either first before the second or the second before the first.
\vskip 4mm
\noindent
{\bf $\cps$ operations}\\
Now, we can define $\cps$-transforms. For $f: X\times Y \ra Z$, we have
\[\cpsl(f)=  \cC(f) \circ \pul_{X, Y}: \cC(X)\times \cC(Y) \lra \cC(Z) \]
given, for $M\in \cC(X)$ and $N\in \cC(Y)$, by
\[ \cpsl(f)(M,N)=\la h_{:\cP(Z)}. M(\la x_{:X}. N(\la y_{:Y}. h(f(x,y)))).\]
Right version is similar.\\
The most popular $\cps$-transforms are those for the evaluation (application) $ ev : X \times (X \Ra Y) \rightarrow Y$
\[\cpsl(ev)=  \cC(ev) \circ \pul_{X,X\Ra Y}: \cC(X)\times \cC(X\Ra Y) \lra \cC(Y)  \]
given, for $M\in \cC(X)$ and $N\in \cC(X\Ra Y)$, by
\[ \cpsl(ev)(M,N)=\la h_{:\cP(Y)}. M(\la x_{:X}. N(\la g_{:X\Ra Y}. h(g\, x))).\]
Right version is similar.
\vskip 4mm
\noindent
{\bf Some functions having products as their domains}\\
There are also other morphisms having useful transforms. Below we list some to introduce notation.

Left evaluations
\[\epsl_X= \la h_{:\cP(X)}.\la x_{:X}. h(x)  : \cP(X) \times X \rightarrow \bt;\]
\[\epsl^{,X}_Y=\epsl_Y= \la c_{:\cP(X\times Y)}.\la y_{:Y}. \la x_{:X}. c(x,y)  : \cP(X\times Y)\times Y \rightarrow \cP(X); \]
and right evaluations
\[ \epsr_X= \la x_{:X}.\la h_{:\cP(X)}. h(x)  : X \times \cP(X) \rightarrow \bt; \]
\[\epsr^{,X}_Y=\epsr_Y=\la y_{:Y}. \la c_{:\cP(X\times Y)}. \la x_{:X}. c(x,y)  : Y\times \cP(X\times Y) \rightarrow \cP(X). \]
$\mos$-operations are the algebraic counterpart of the familiar interpretation of generalized quantifiers of Mostowski (again, we give definitions for total and partial case). Left $\mos$'es
\[ \mosl_X =\la Q_{:\cC(X)}.\la c_{:\cP(X)}. Q(c): \cC(X)\times \cP(X) \rightarrow \textbf{t}; \]
\[ \mosl_Y =\la Q_{:\cC(Y)}.\la c_{:\cP(X\times Y)}. \la x_{:X}. Q( \la y_{:Y}. c(x,y)) : \cC(Y)\times \cP(X\times Y) \rightarrow \cP(X);\]
and right $\mos$'es
\[ \mosr_X =\la c_{:\cP(X)}. \la Q_{:\cC(X)}. Q(c) :\cP(X) \times \cC(X) \rightarrow \textbf{t}; \]
\[ \mosr_Y =\la c_{:\cP(X\times Y)}. \la Q_{:\cC(Y)}. \la x_{:X}. Q( \la y_{:Y}. c(x,y)) :\cP(X\times Y) \times \cC(Y) \rightarrow \cP(X).\]

\section{Continuation-Based In Situ Strategies}

Below we illustrate Strategy \textbf{C}, \textbf{D} $\&$ \textbf{E} on examples involving one, two and three QPs (for the description of Strategies \textbf{A} $\&$ \textbf{B}, see \cite{GZ}). In each strategy, the leaves in the Computation Trees have the same labels: QPs are interpreted as $\cC$-computations, and predicates are interpreted as lifted (`continuized') relations. The main difference among the three approaches consists in the operations ($\cps$'es, $\sh$'es, $\eps$'es) used as labels of the inner nodes of the Computation Trees.

\subsection{Strategy C}

\textbf{Sentence with one QP}, e.g. \textit{Every kid (most kids) entered.}

\begin{flushleft}
(\textbf{C1}) Surface Structure Tree and the corresponding Computation Tree

\medskip
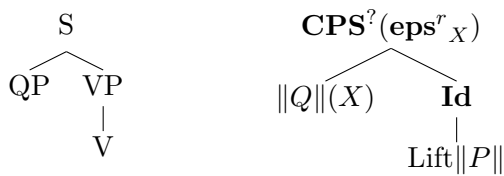
\Tree [.S [.QP ][.VP V ]]
\Tree [.$\cps^?(\epsr_{X})$ [.$\|Q\|(X)$  ][.$\Id$ Lift$\|P\|$ ]]
\par
\medskip
\end{flushleft}
with $\Id$ standing for the identity operation. We use $\cps^?$ when it does not matter whether we apply $\cpsl$ or $\cpsr$. This is the case when one of the arguments is a lifted element (like interpretations of predicates in this strategy). \textbf{Strategy C} yields one reading for a sentence with one QPs.

\vskip 2mm

\noindent \textbf{Sentence with two QPs}, e.g. \textit{Every girl likes a boy.}

\begin{flushleft}
(\textbf{C2}) Surface Structure Tree and the corresponding Computation Tree

\medskip
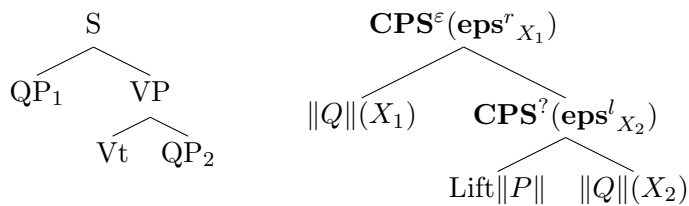
\Tree [.S [.QP_1 ][.VP [.Vt ] [.QP_2 ]]]
\Tree [.$\cps^\varepsilon(\epsr_{X_1})$ [.$\|Q\|(X_1)$ ][.$\cps^?(\epsl_{X_2})$ [.Lift$\|P\|$ ][.$\|Q\|(X_2)$ ]]]
\par
\end{flushleft}
with $\varepsilon\in \{l,r\}$. Depending on whether we use $\cpsl$ or $\cpsr$, we get either one or the other of the two asymmetric readings for a sentence with two QPs. Thus \textbf{Strategy C} yields two readings for a sentence with two QPs, corresponding to the two $\cps$'es.

\vskip 2mm

\noindent \textbf{Sentence with three QPs}, e.g. \textit{Some teacher gave every student most books.}

\begin{flushleft}
(\textbf{C3}) Surface Structure Tree and the corresponding Computation Tree

\medskip
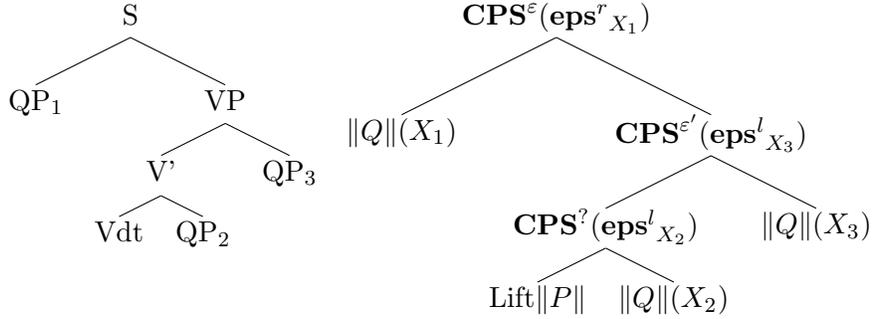
\Tree [.S [.QP_1 ][.VP [.V' [.Vdt ][.QP_2 ]][.QP_3 ]]]
\Tree [.$\cps^\varepsilon(\epsr_{X_1})$ [.$\|Q\|(X_1)$ ][.$\cps^{\varepsilon'}(\epsl_{X_3})$ [.$\cps^?(\epsl_{X_2})$ [.Lift$\|P\|$ ][.$\|Q\|(X_2)$ ]][.$\|Q\|(X_3)$ ]]]
\par
\end{flushleft}
\textbf{Strategy C} provides four asymmetric readings for a sentence with a ditransitive verb such that QP in subject position can be placed either first or last only, corresponding to the four possible combinations of the two $\cps$'es.

\textbf{Strategy C} allows a uniform in situ analysis of quantifiers. However, it cannot be straightforwardly extended to account for sentences involving 3 QPs - as discussed in \cite{bekki} and proved in \cite{GZ}, it only provides four out of six readings available for such sentences. Below we define \textbf{Strategy D}, a minimally augmented empirically adequate version of \textbf{Strategy C}.

\subsection{Strategy D}

To include the two readings missing from \textbf{Strategy C}, \textbf{Strategy D} combines $\cps$'es to define two new operations: $\shl$ and $\shr$ (one can also combine $\pu$- and $\cps$-operations to define those new operations, this is left for another place though).

\vskip 2mm

\noindent \textbf{Sentence with one QP}, e.g. \textit{Every kid (most kids) entered.}

\begin{flushleft}
(\textbf{D1}) Surface Structure Tree and the corresponding Computation Tree

\medskip
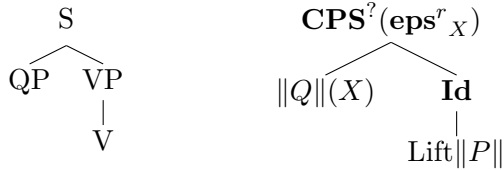
\Tree [.S [.QP ][.VP V ]]
\Tree [.$\cps^?(\epsr_{X})$ [.$\|Q\|(X)$  ][.$\Id$ Lift$\|P\|$ ]]
\par
\medskip
\end{flushleft}
Just as in \textbf{Strategy C}, \textbf{Strategy D} yields one reading for a sentence with one QPs.

\vskip 2mm

\noindent \textbf{Sentence with two QPs}, e.g. \textit{Every girl likes a boy.}

\begin{flushleft}
(\textbf{D2}) Surface Structure Tree and the corresponding Computation Tree

\medskip
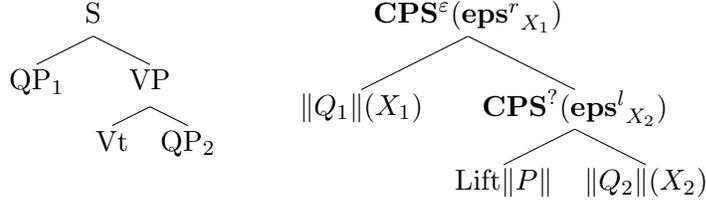
\Tree [.S [.QP_1 ][.VP [.Vt ] [.QP_2 ]]]
\Tree [.$\cps^\varepsilon(\epsr_{X_1})$ [.$\|Q_1\|(X_1)$ ][.$\cps^?(\epsl_{X_2})$ [.Lift$\|P\|$ ][.$\|Q_2\|(X_2)$ ]]]
\par
\end{flushleft}
Just as in \textbf{Strategy C}, \textbf{Strategy D} yields both asymmetric readings for such sentences.

\vskip 2mm

\noindent \textbf{Sentence with three QPs}, e.g. \textit{Some teacher gave every student most books.}

\begin{flushleft}
(\textbf{D3}) Surface Structure Tree and the corresponding Computation Tree

\medskip
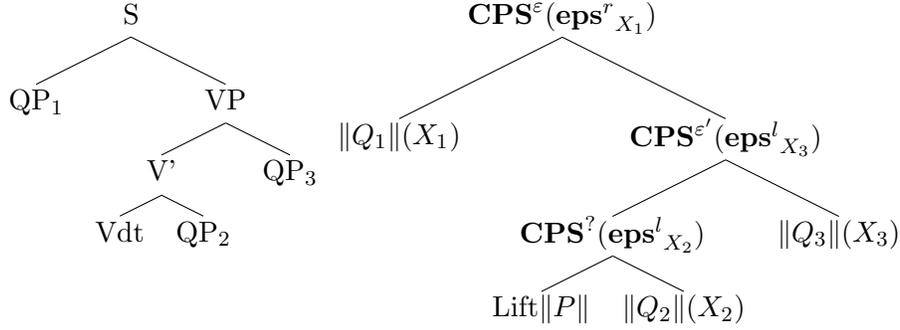
\Tree [.S [.QP_1 ][.VP [.V' [.Vdt ][.QP_2 ]][.QP_3 ]]]
\Tree [.$\cps^\varepsilon(\epsr_{X_1})$ [.$\|Q_1\|(X_1)$ ][.$\cps^{\varepsilon'}(\epsl_{X_3})$ [.$\cps^?(\epsl_{X_2})$ [.Lift$\|P\|$ ][.$\|Q_2\|(X_2)$ ]][.$\|Q_3\|(X_3)$ ]]]
\par
\end{flushleft}
Just as in \textbf{Strategy C}, the Computation Tree above gives rise to the four asymmetric readings for sentences with ditransitive verbs such that QP in subject position can be placed either first or last only. To get the first missing reading
\[QP_3 > QP_1 > QP_2, \hskip 10cm \]
we define a new operation
\[ \shr: \cC \cP (X_1 \times X_3) \times \cC (X_3) \rightarrow \cP \cC (X_1) \]
such that
\[ \shr(S_2,S_3) = \lambda {S_1}_{: \cC (X_1)}. \cpsl(\epsr_{X_3}) (S_3, \cpsl(\epsr_{X_1})(S_1, S_2))(id_t) \hskip 10cm\]
for $S_2\in \cC \cP (X_1 \times X_3)$  and $S_3\in \cC (X_3)$.
\vskip 2mm
\noindent The corresponding computation tree is now as follows
\begin{flushleft}
(\textbf{D3'}) Computation Tree

\medskip
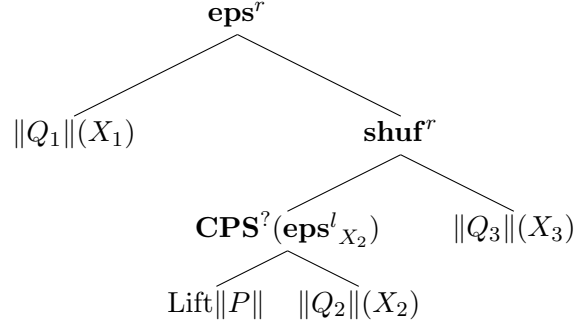
\Tree [.$\epsr$ [.$\|Q_1\|(X_1)$ ][.$\shr$ [.$\cps^?(\epsl_{X_2})$ [.Lift$\|P\|$ ][.$\|Q_2\|(X_2)$ ]][.$\|Q_3\|(X_3)$ ]]]
\par
\end{flushleft}
To get the second of the two missing readings
\[QP_2 > QP_1 > QP_3, \hskip 10cm \]
we define a new operation
\[ \shl : \cC \cP (X_1 \times X_3) \times \cC (X_3) \rightarrow \cP \cC (X_1) \]
\[ \shl(S_2,S_3) = \lambda {S_1}_{: \cC (X_1)}.\cpsl(\epsl_{X_3}) (\cpsl(\epsl_{X_1}) (S_2, S_1), S_3)(id_t) \hskip 10cm\]
for  $S_2\in \cC \cP (X_1 \times X_3)$ and $S_3\in \cC (X_3)$.
\vskip 2mm
\noindent The corresponding computation tree is now as follows
\begin{flushleft}
(\textbf{D3''}) Computation Tree

\medskip
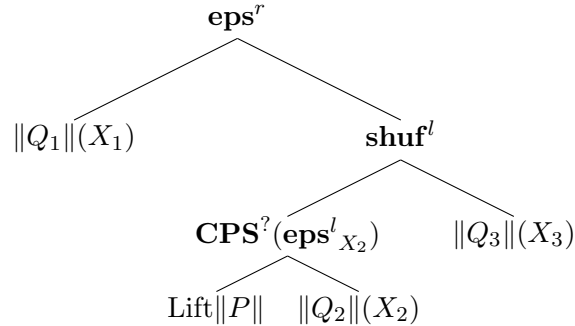
\Tree [.$\epsr$ [.$\|Q_1\|(X_1)$ ][.$\shl$ [.$\cps^?(\epsl_{X_2})$ [.Lift$\|P\|$ ][.$\|Q_2\|(X_2)$ ]][.$\|Q_3\|(X_3)$ ]]]
\par

\end{flushleft}
In \textbf{Strategy D}, unlike in \textbf{Strategy C}, we get all the asymmetric readings. Thus \textbf{Strategy D} is both in situ and free of the empirical deficiencies in \textbf{Strategy C}. This is achieved at the price, though, of extending the list of operations (adopted in strategy C) by two more involved operations: left and right $\sh$. Moreover, the process translating Surface Structure Trees into the Semantic (Computation) Trees loses its uniformity - depending on the class of sentences (involving one, two or three QPs) and the particular readings considered, different semantic operations are used. Below we define \textbf{Strategy E}, a maximally uniformized version of \textbf{Strategy D}.

\subsection{Strategy E}

In \textbf{Strategy E}, for each asymmetric reading its respective $\sh$-operation will be defined. We will index $\sh$'es with permutations of the $n$ distinct QPs involved, i.e., $\sh^\sigma$ will denote the operation that determines the reading of the sentence with the QPs ordered according to the permutation $\sigma$.

\noindent \textbf{Sentence with one QP}, e.g. \textit{Every kid (most kids) entered.}

\begin{flushleft}
(\textbf{E1}) Surface Structure Tree and the corresponding Computation Tree

\medskip
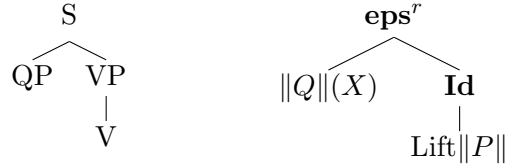
\Tree [.S [.QP ][.VP V ]]
\Tree [.$\epsr$ [.$\|Q\|(X)$  ][.$\Id$ Lift$\|P\|$ ]]
\par
\medskip
\end{flushleft}
\textbf{Strategy E} yields one reading for a sentence with one QPs.

\vskip 2mm

\noindent \textbf{Sentence with two QPs}, e.g. \textit{Every girl likes a boy.}

\begin{flushleft}
(\textbf{E2}) Surface Structure Tree

\medskip
\Tree [.S [.QP_1 ][.VP [.Vt ] [.QP_2 ]]]

\par
\end{flushleft}

\[QP_1 > QP_2 \hskip 14cm\]
\[ \sh^{1,2} : \cC \cP (X_1 \times X_2) \times \cC (X_2) \rightarrow \cP \cC (X_1) \]
\[ \sh^{1,2}(S, S_2) = \lambda {S_1}_{: \cC (X_1)}. \cpsl(\epsr_{X_1}) (S_1, \cpsl(\epsl_{X_2})(S, S_2))(id_t) \hskip 10cm\]
for $S \in \cC \cP (X_1 \times X_2)$  and $S_2\in \cC (X_2)$.
\[QP_2 > QP_1 \hskip 14cm\]
\[ \sh^{2,1} : \cC \cP (X_1 \times X_2) \times \cC (X_2) \rightarrow \cP \cC (X_1) \]
\[ \sh^{2,1}(S, S_2) = \lambda {S_1}_{: \cC (X_1)}. \cpsl(\epsl_{X_1}) (\cpsl(\epsl_{X_2})(S, S_2), S_1)(id_t) \hskip 10cm\]
for $S \in \cC \cP (X_1 \times X_2)$  and $S_2\in \cC (X_2)$.
\vskip 2mm
\noindent The corresponding computation trees are now as follows
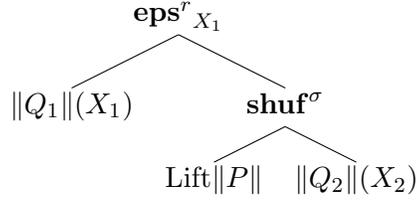
\begin{flushleft}

\Tree [.$\epsr_{X_1}$ [.$\|Q_1\|(X_1)$ ][.$\sh^\sigma$ [.Lift$\|P\|$ ][.$\|Q_2\|(X_2)$ ]]]

\medskip

\end{flushleft}
where $\sigma$ is a permutation of $\{1, 2\}$.
\vskip 2mm
\noindent \textbf{Sentence with three QPs}, e.g. \textit{Some teacher gave every student most books.}

\begin{flushleft}
(\textbf{E3}) Surface Structure Tree

\medskip
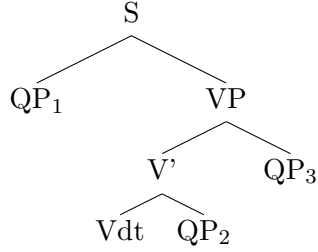
\Tree [.S [.QP_1 ][.VP [.V' [.Vdt ][.QP_2 ]][.QP_3 ]]]
\par
\end{flushleft}
For each reading, we now define its respective $\sh$-operation
\[QP_1 > QP_2 > QP_3 \hskip 10cm\]
\[ \sh^{1,2,3} : \cC \cP (X_1 \times X_3) \times \cC (X_3) \rightarrow \cP \cC (X_1) \]
\[ \sh^{1,2,3}(S_2,S_3) = \lambda {S_1}_{: \cC (X_1)}. \cpsl(\epsr_{X_1}) (S_1, \cpsl(\epsl_{X_3})(S_2, S_3))(id_t)  \]
for $S_2\in \cC \cP (X_1 \times X_3)$  and $S_3\in \cC (X_3)$.
\[QP_3 > QP_2 > QP_1 \hskip 10cm \]
\[ \sh^{3,2,1} : \cC \cP (X_1 \times X_3) \times \cC (X_3) \rightarrow \cP \cC (X_1)  \]
\[\sh^{3,2,1}(S_2,S_3) = \lambda {S_1}_{: \cC (X_1)}.\cpsl(\epsl_{X_1}) (\cpsl(\epsr_{X_3}) (S_3, S_2), S_1)(id_t) \]
for  $S_2\in \cC \cP (X_1 \times X_3)$ and $S_3\in \cC (X_3)$.
\[QP_1 > QP_3 > QP_2  \hskip 10cm \]
\[ \sh^{1,3,2} : \cC \cP (X_1 \times X_3) \times \cC (X_3) \rightarrow \cP \cC (X_1)  \]
\[ \sh^{1,3,2}(S_2,S_3) = \lambda {S_1}_{: \cC (X_1)}. \cpsl(\epsr_{X_1}) (S_1, \cpsl(\epsr_{X_3})(S_3, S_2))(id_t) \]
for $S_2\in \cC \cP (X_1 \times X_3)$  and $S_3\in \cC (X_3)$.
\[QP_2 > QP_3 > QP_1  \hskip 10cm \]
\[ \sh^{2,3,1} : \cC \cP (X_1 \times X_3) \times \cC (X_3) \rightarrow \cP \cC (X_1)  \]
\[ \sh^{2,3,1}(S_2,S_3) = \lambda {S_1}_{: \cC (X_1)}.\cpsl(\epsl_{X_1}) (\cpsl(\epsl_{X_3}) (S_2, S_3), S_1)(id_t) \]
for  $S_2\in \cC \cP (X_1 \times X_3)$ and $S_3\in \cC (X_3)$.
\[QP_3 > QP_1 > QP_2 \hskip 10cm \]
\[ \sh^{3,1,2} : \cC \cP (X_1 \times X_3) \times \cC (X_3) \rightarrow \cP \cC (X_1)  \]
\[ \sh^{3,1,2}(S_2,S_3) = \lambda {S_1}_{: \cC (X_1)}. \cpsl(\epsr_{X_3}) (S_3, \cpsl(\epsr_{X_1})(S_1, S_2))(id_t) \]
for $S_2\in \cC \cP (X_1 \times X_3)$  and $S_3\in \cC (X_3)$.
\[QP_2 > QP_1 > QP_3  \hskip 10cm  \]
\[ \sh^{2,1,3} : \cC \cP (X_1 \times X_3) \times \cC (X_3) \rightarrow \cP \cC (X_1)  \]
\[ \sh^{2,1,3}(S_2,S_3) = \lambda {S_1}_{: \cC (X_1)}.\cpsl(\epsl_{X_3}) (\cpsl(\epsl_{X_1}) (S_2, S_1), S_3)(id_t) \]
for  $S_2\in \cC \cP (X_1 \times X_3)$ and $S_3\in \cC (X_3)$.
\vskip 2mm
\noindent The corresponding computation trees are now as follows
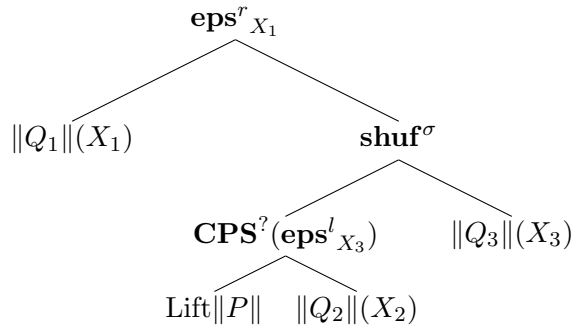
\begin{flushleft}
\medskip
\Tree [.$\epsr_{X_1}$ [.$\|Q_1\|(X_1)$ ][.$\sh^\sigma$ [.$\cps^?(\epsl_{X_3})$ [.Lift$\|P\|$ ][.$\|Q_2\|(X_2)$ ]][.$\|Q_3\|(X_3)$ ]]]
\par
\end{flushleft}
where $\sigma$ is a permutation of $\{1, 2, 3\}$.\\

\noindent \textbf{Strategy E} is both in situ and uniform, i.e., regardless of the class of sentences and the particular readings considered, it uses the same semantic operations ($\cps$'es, $\sh$'es and $\eps$'es). Compared to \textbf{Strategy D}, however, it fares worse with respect to the requirement (5) that the interpretation process operate `on the fly', i.e., the particular reading(s) of a sentence should be determined as late as possible. In \textbf{Strategy D}, it is only the two readings missing from \textbf{Strategy C} that are predetermined by the $\sh$'es. In \textbf{Strategy E}, it is all of the readings for sentences involving 2 or 3 QPs that are determined by the $\sh$'es. So between \textbf{Strategy D} and \textbf{Strategy E} there is a trade-off of gains and costs: having `on the fly' process vs. uniformity in semantic operations. Perhaps, if it could be empirically shown that some readings take longer to process (e.g., the two readings missing from \textbf{Strategy C}), then one could hypothesize that they also involve comparably more difficult semantic operations (i.e., $\sh$'es, in this case). Such empirical findings could be then taken to support \textbf{Strategy D} over \textbf{E}. Obviously, the empirical findings could be also found to support some mixed-strategy, located somewhere between \textbf{Strategy D} and \textbf{Strategy E}. One should also notice that in the in situ strategies discussed in this paper (whether or not the readings are predetermined) the arguments of the semantic operations applied are used `on the spot', i.e., unlike in Cooper's Storage mechanism (\cite{cooper}), the arguments do not get stored and retrieved when needed.

\section{Conclusion}

Recent continuation-based semantics provide a compelling approach to quantification for providing a non-movement (in situ) analysis of quantifiers. In this paper, we have discussed three such possible continuation-based strategies, together with their relative merits and costs.

\section{Appendix: Strong Monads}

A {\em monad} on $Set$ is a triple $(T,\eta,\mu)$, where $T: Set \lra Set$ is an endofunctor, $\eta : 1_{Set} \lra T$ and $\mu: T^2 \lra T$ are natural transformations making the following diagrams
\begin{center}
\xext=1600 \yext=800 \adjust[`I;I`;I`;`I]
\begin{picture}(\xext,\yext)(\xoff,\yoff)
\settriparms[-1`1`1;400]
\put(480,190){$\mu$}
\putAtriangle(0,0)[T^2`T`T;\eta_T``1_T ]

\settriparms[1`-1`-1;400]
\put(1060,190){$\mu$}
\putAtriangle(800,0)[T^2`\phantom{T}`T;`T(\eta)`1_T ]

\settriparms[1`1`0;400]
\putAtriangle(400,400)[T^3`\phantom{T^2}`\phantom{T^2};\mu_T`T(\mu)` ]
\end{picture}
\end{center}
%
%
%
%
%
%
%
commute. $\eta$ and $\mu$ are {\em unit} and {\em multiplication} of the monad $T$, respectively.

In order to have a well behaved notion of computation, a monad has to be strong, c.f. \cite{kock1}, \cite{kock4}, \cite{moggi}. Fortunately, all monads on $Set$ are strong. Technically, we need the monad to be bi-strong as we will need to `extend computations' both from the left and from the right. As the binary product (the only tensor that we consider in $Set$) is commutative, any strong monad in $Set$ is bi-strong.

Let $(T,\eta,\mu)$ be a monad on $Set$. The {\em left strength} on  $(T,\eta,\mu)$ is a natural transformation with components
\[  \stl_{X,Y} : T(X)\times Y \lra T(X\times Y) \]
for sets $X$ and $Y$, making the diagrams
\begin{center}
\xext=1100 \yext=600 \adjust[`I;I`;I`;`I]
\begin{picture}(\xext,\yext)(\xoff,\yoff)
\settriparms[1`1`-1;550]
  \putVtriangle(0,0)[T(X)\times Y\times Z`T(X\times Y\times Z)`T(X\times Y)\times Z;\stl_{X,Y\times Z}`\stl_{X,Y}\times 1` \stl_{X\times Y, Z}]
\end{picture}
\end{center}
and
\begin{center}
\xext=2400 \yext=1100
\begin{picture}(\xext,\yext)(\xoff,\yoff)
\putmorphism(0,1000)(0,-1)[X\times Y`\phantom{T(X)\times Y}`\eta_X \times 1]{500}{1}l

\putmorphism(50,1000)(2,-1)[\phantom{X\times Y}`\phantom{T(X\times Y)}`\eta_{X\times Y}]{1000}{1}r

 \setsqparms[1`-1`0`1;1200`500]
\putsquare(0,0)[T(X)\times Y`T(X\times Y)`T^2(X)\times Y`T(T(X)\times Y);{\stl}_{X,Y}`\mu_X\times 1``\stl_{T(X),Y}]

\putmorphism(1250,500)(2,-1)[\phantom{T(X\times Y)}`\phantom{T^2(X\times Y)}`\mu_{X\times Y}]{1000}{-1}r
\putmorphism(1200,0)(1,0)[\phantom{T(X\times T(Y))}`T^2(X\times Y)`T({\stl}_{X\times Y})]{1200}{1}b

\end{picture}
\end{center}
commute.

The {\em right strength} is a natural transformation with components
\[  \str_{X,Y} : X\times T(Y) \lra T(X\times Y) \]
for sets $X$ and $Y$, making the diagrams
\begin{center}
\xext=1100 \yext=600 \adjust[`I;I`;I`;`I]
\begin{picture}(\xext,\yext)(\xoff,\yoff)
\settriparms[1`1`-1;550]
  \putVtriangle(0,0)[X\times Y\times T(Z)`T(X\times Y\times Z)`X\times T(Y\times Z);\str_{X\times Y,Z}`1\times \str_{Y,Z}` \str_{X,Y\times Z}]
\end{picture}
\end{center}
and
\begin{center}
\xext=2400 \yext=1100
\begin{picture}(\xext,\yext)(\xoff,\yoff)
\putmorphism(0,1000)(0,-1)[X\times Y`\phantom{X\times T(Y)}`1\times \eta_Y]{500}{1}l

\putmorphism(50,1000)(2,-1)[\phantom{X\times Y}`\phantom{T(X\times Y)}`\eta_{X\times Y}]{1000}{1}r

 \setsqparms[1`-1`0`1;1200`500]
\putsquare(0,0)[X\times T(Y)`T(X\times Y)`X\times T^2(Y))`T(X\times T(Y));{\str}_{X,Y}`1\times \mu_Y``\str_{X,T(Y)}]

\putmorphism(1250,500)(2,-1)[\phantom{T(X\times Y)}`\phantom{T^2(X\times Y)}`\mu_{X\times Y}]{1000}{-1}r
\putmorphism(1200,0)(1,0)[\phantom{T(X\times T(Y))}`T^2(X\times Y)`T({\str}_{X\times Y})]{1200}{1}b

\end{picture}
\end{center}
commute.

The monad $(T,\eta,\mu)$ on $Set$ together with two natural transformations $\stl$ and $\str$ of right and left strength is a {\em bi-strong monad} if, for any sets $X$, $Y$, $Z$, the square
\begin{center}
\xext=1200 \yext=600 \adjust[`I;I`;I`;`I]
\begin{picture}(\xext,\yext)(\xoff,\yoff)
\setsqparms[1`1`1`1;1200`500]
\putsquare(00,0)[X\times T(Y)\times Z`X\times T((Y\times Z)`T(X\times Y)\times Z`T(X\times Y \times Z);1_X\times \stl_{Y,Z}`\str_{X,Y}\times 1_Z`\str_{X,Y\times Z}`\stl_{X\times Y, Z}]
\end{picture}
\end{center}
commutes.

\vskip 2mm
{\em Remarks.}
\begin{enumerate}
  \item The general definition for a strong functor on a monoidal category $(\cC,\otimes,I,\alpha, \lambda, \rho)$ contains a yet another diagram concerning the unit. As the only tensor considered is the cartesian product, this part of the general definition is irrelevant.
  \item If a monad on a symmetric monoidal category is strong, then (having, say, left strength only) one can - using symmetry - easily define the right strength making it a bi-strong monad. In that sense, the concept of a bi-strong monad is redundant in the symmetric monoidal categories. However, since the order does matter in the natural language, we prefer to make both strengths explicitly given as a part of the structure.
\end{enumerate}
\vskip 2mm

As  already mentioned, each monad $(T,\eta,\mu)$ on $Set$ is bi-strong. We shall define the right and left strength. Fix sets $X$ and $Y$.
For $x\in X$ and $y\in Y$, we have functions
\[ l_{y} : X\lra X\times Y, \hskip 5mm{\rm and} \hskip 5mm  r_{x} : Y\lra X\times Y,\]
such that
\[ l_{y}(x) = \lk x,y\rk,  \hskip 5mm{\rm and} \hskip 5mm r_{x}(y) = \lk x,y\rk. \]
The left and right strength
\[ \stl_{X,Y}: T(X)\times Y \lra T(X\times Y) \hskip 5mm{\rm and} \hskip 5mm \str_{X,Y}: X\times T(Y) \lra T(X\times Y) \]
are given for $x\in X$, $s\in T(X)$, $y\in Y$ and $t\in T(Y)$ by
\[ \stl_{X,Y}(s,y)= T(l_y)(s)   \hskip 5mm{\rm and} \hskip 5mm  \str_{X,Y}(x,t)= T(r_x)(t), \]
respectively. We drop indices $_{X,Y}$ when it does not lead to confusion.

One can check that the monad $(T,\eta,\mu)$ equipped with so defined natural transformations $\stl$ and  $\str$ is bi-strong.


\begin{thebibliography}{}


\bibitem {barker}
Barker, C.:
Continuations and the nature of quantification.
Natural Language Semantics 10, 211-242 (2002).

\bibitem {barker:shan}
Barker, C., Shan, C.c.:
Continuations and Natural Language.
Oxford University Press (2014).

\bibitem {bar:coop}
Barwise, J., Cooper, R.:
Generalized Quantifiers and Natural Language.
Linguistics \& Philosophy 4, 159-219 (1981).

\bibitem {bekki}
Bekki, D., Assai, K.:
Representing Covert Movements by Delimited Continuations
In Nakakoji, K., Murakami, Y., McCready, E. (eds.) New Frontiers in Artificial Intelligence, JSAI-isAI, LNAI 6284, 161-180, (2010).

\bibitem {benthem}
Benthem, J.:
Polyadic quantifiers.
Linguistics \& Philosophy 12, 437-464 (1989).

\bibitem {cooper}
Cooper, R.:
Quantification and Syntactic Structure.
Dordrecht: Reidel (1983).

\bibitem {groote}
de Groote, P.:
Type raising, continuations, and classical logic.
In van Rooy, R., Stokhof, M. (eds.), Proceedings of the 13th Amsterdam Colloquium, Institute
for Logic, Language and Computation, Universiteit van Amsterdam 97-101 (2001).


\bibitem {keenan87}
Keenan, E. L.:
Unreducible n-ary quantifiers in natural language.
In Gärdenfors P. (ed.), Generalized Quantifier: Linguistic and Logical Approaches, Reidel Dordrecht. 109-150 (1987).


\bibitem {kiselyov:shan}
Kiselyov, O., Shan, C.c.:
Continuation Hierarchy and Quantifier Scope.
In McCready, E., Yabushita, K., Yoshimoto, K. (eds.), Formal Approaches to Semantics and Pragmatics: Japanese and Beyond.
Studies in Linguistics and Philosophy, Springer Netherlands, 105-134 (2014).

\bibitem {kock1}
Kock, A.:
Monads on symmetric monoidal closed categories.
Arch. Math. (Basel), 21:1–10 (1970).

\bibitem {kock4}
Kock, A.:
Strong functors and monoidal monads.
Arch. Math. (Basel), 23:113–120 (1972).

\bibitem {lind}
Lindstr\"{o}m, P.:
First-order predicate logic with generalized quantifiers.
Theoria 32, 186-95.(1966).

\bibitem {may}
May, R.:
The Grammar of Quantification.
PhD dissertation, MIT (1977).

\bibitem {may85}
May, R.:
Logical Form: Its Structure and Derivation.
MIT Press (1985).

\bibitem {moggi}
Moggi, E.:
The notion of computation and monads.
Information And Computation 93:55–92 (1991).

\bibitem {montague}
Montague, R.:
The proper treatment of quantification in ordinary English.
In: Thomason, R. (ed.), Formal Philosophy: Selected Papers of Richard Montague.
New Haven and London: Yale University Press, 247-271 (1974).

\bibitem {most}
Mostowski, A.:
On a generalization of quantifiers.
Fundamenta Mathematicae 44, 12-36 (1957).

\bibitem {szab}
Szabolcsi, A.:
Quantification.
Cambridge University Press, Cambridge (2010).

\bibitem {BZ}
Zawadowski, M.:
Formalization of the feature system in terms of pre-orders.
In Bellert, I, Feature System for Quantification Structures in Natural Language, Foris Dordrecht, 155-175 (1989).

\bibitem {GZ}
Grudzinska, J; Zawadowski, M.:
Scope ambiguities, monads and strengths.
arXiv:1605.03981.

\end{thebibliography}
\end{document}